\documentclass{article}
\usepackage{amssymb}
\usepackage{amsfonts}
\usepackage{amsmath}
\usepackage{epsfig}

\begin{document}

\begin{center}
{\Large About some properties of algebras obtained by the Cayley-Dickson
process}%
\begin{equation*}
\end{equation*}

Cristina FLAUT%
\begin{equation*}
\end{equation*}
\end{center}

\bigskip

\textbf{Abstract.} {\small This paper is a short survey} {\small about some
properties of algebras obtained by the Cayley-Dickson process and some of
their applications.}

\begin{equation*}
\end{equation*}

\begin{equation*}
\end{equation*}

\begin{equation*}
\end{equation*}%
\begin{equation*}
\end{equation*}

\textbf{1.} \textbf{Introduction}%
\begin{equation*}
\end{equation*}%
\begin{equation*}
\end{equation*}

It is well known that in October 1843, William Rowan Hamilton made a great
discovery finding quaternion algebra, a $4$-dimensional algebra over $%
\mathbb{R}$ which is an associative and \ a noncommutative algebra. In
December $1843,$ \ John Graves discovered the octonions, an $8$-dimensional
algebra over $\mathbb{R}$ which is nonassociative and noncommutative
algebra. These algebras were later rediscovered \ by Arthur Cayley in $1845$
and are also known sometimes as the \textit{Cayley numbers.} This process,
of passing from $\mathbb{R}$ to $\mathbb{C}$, from $\mathbb{C}$ to $\mathbb{H%
}$ and from$\ \mathbb{H}$ to $\mathbb{O}$ has been generalized to algebras
over fields and over rings. It is called the \ \textit{Cayley-Dickson
doubling process} \ or the \textit{Cayley--Dickson process.}

Even if are old, Quaternion and Octonion algebras have at \ present many
applications, especially in physics, coding theory, computer science, etc.
For example, reliable high rate of transmission can be obtained using
Space-Time coding. For constructing Space-Time codes, Quaternion division
algebras were chosen as a new tool, as for example the Alamouti code, which
can be built from a quaternion division algebra (see [Al; 98]).

The classical Cayley-Dickson process is briefly presented in the following.
For details about this, the reader is referred to $\left[ \text{Sc; 66}%
\right] $ and [Sc; 54]. From now on, in the whole paper, we will consider $K$
a field with $charK\neq 2.$

Let $\ A$ be an algebra over the field $K.$ A unitary algebra $A\neq K$ such
that we have $x^{2}+\alpha \left( x\right) x+\beta \left( x\right) =0,$ for
each $x\in A,$ with $\alpha \left( x\right) ,\beta \left( x\right) \in K,$
is called a \textit{quadratic algebra}.

Let $A$ be a finite dimensional unitary algebra over a field $\ K$ with a 
\textit{scalar} \textit{involution} $\,$%
\begin{equation*}
\,\,\,\overline{\phantom{x}}:A\rightarrow A,a\rightarrow \overline{a},
\end{equation*}%
$\,\,$ i.e. a linear map satisfying the following relations:$\,\,\,\,\,$%
\begin{equation*}
\overline{ab}=\overline{b}\overline{a},\,\overline{\overline{a}}=a,
\end{equation*}%
$\,\,$and 
\begin{equation*}
a+\overline{a},a\overline{a}\in K\cdot 1\ \text{for all }a,b\in A.\text{ }
\end{equation*}%
An element $\,\overline{a}$ is called the \textit{conjugate} of the element $%
a,$ the linear form$\,\,$%
\begin{equation*}
\,\,t:A\rightarrow K\,,\,\,t\left( a\right) =a+\overline{a}
\end{equation*}%
and the quadratic form 
\begin{equation*}
n:A\rightarrow K,\,\,n\left( a\right) =a\overline{a}\ 
\end{equation*}%
are called the \textit{trace} and the \textit{norm \ }of \ the element $a,$
respectively$.$ Therefore, such an algebra $A$ with a scalar involution is
quadratic. $\,$

Let$\,\,\,\gamma \in K$ \thinspace be a fixed non-zero element. On the
vector space $A\oplus A$, we define the following algebra multiplication: 
\begin{equation}
A\oplus A:\left( a_{1},a_{2}\right) \left( b_{1},b_{2}\right) =\left(
a_{1}b_{1}+\gamma \overline{b_{2}}a_{2},a_{2}\overline{b_{1}}%
+b_{2}a_{1}\right) .  \tag{1.1}
\end{equation}%
\newline
We obtain an algebra structure over $A\oplus A,$ denoted by $\left( A,\gamma
\right) $ and called the \textit{algebra obtained from }$A$\textit{\ by the
Cayley-Dickson process.} $\,$It results that $\dim \left( A,\gamma \right)
=2\dim A$.

For $x\in \left( A,\gamma \right) $, $x=\left( a_{1},a_{2}\right) $, the map 
\begin{equation}
\,\,\,\overline{\phantom{x}}:\left( A,\gamma \right) \rightarrow \left(
A,\gamma \right) \,,\,\,x\rightarrow \bar{x}\,=\left( \overline{a}_{1},\text{%
-}a_{2}\right) ,  \tag{1.2}
\end{equation}%
\newline
is a scalar involution of the algebra $\left( A,\gamma \right) $, extending
the involution $\overline{\phantom{x}}\,\,\,$of the algebra $A.$ Let 
\begin{equation*}
\,t\left( x\right) =t(a_{1})
\end{equation*}%
and$\,\,\,$ 
\begin{equation*}
n\left( x\right) =n\left( a_{1}\right) -\gamma n(a_{2})
\end{equation*}%
be $\,\,$the \textit{trace} and the \textit{norm} of the element $x\in $ $%
\left( A,\gamma \right) ,$ respectively.\thinspace $\,$

\thinspace If we take $A=K$ \thinspace and apply this process $t$ times, $%
t\geq 1,\,\,$we obtain an algebra over $K,\,\,$%
\begin{equation}
A_{t}=\left( \frac{\alpha _{1},...,\alpha _{t}}{K}\right) .  \tag{1.3}
\end{equation}%
By induction in this algebra, the set $\{1,e_{2},...,e_{n}\},n=2^{t},$
generates a basis with the properties:%
\begin{equation}
e_{i}^{2}=\alpha _{i}1,\,\,\alpha _{i}\in K,\alpha _{i}\neq 0,\,\,i=2,...,n 
\tag{1.4}
\end{equation}%
and \ 
\begin{equation}
e_{i}e_{j}=-e_{j}e_{i}=\beta _{ij}e_{k},\,\,\beta _{ij}\in K,\,\,\beta
_{ij}\neq 0,i\neq j,i,j=\,\,2,...n,  \tag{1.5}
\end{equation}%
$\ \beta _{ij}$ and $e_{k}$ being uniquely determined by $e_{i}$ and $e_{j}.$

From [Sc; 54], Lemma 4, it results that in any algebra $A_{t}$ with the
basis \newline
$\{1,e_{2},...,e_{n}\}$ satisfying the above relations we have:

\begin{equation}
e_{i}\left( e_{i}x\right) =\alpha _{i}^{2}=(xe_{i})e_{i},  \tag{1.6}
\end{equation}%
for \ every $x\in A$ and for all $i\in \{1,2,...,n\}.$

A finite-dimensional algebra $A$ is a division algebra if and only if $A$
does not contain zero divisors (see [Sc;66]).

An algebra $A$ is called \textit{central simple} if \ the algebra\ $%
A_{F}=F\otimes _{K}A$ is simple for every extension $F$ of $K.$ An algebra $%
A $ is called \textit{alternative} if $x^{2}y=x\left( xy\right) $ and $%
xy^{2}=\left( xy\right) y,$ for all $x,y\in A.$ An algebra $A$ is called 
\textit{\ flexible} if $x\left( yx\right) =\left( xy\right) x=xyx,$ for all $%
x,y\in A$ and \textit{power associative} if the subalgebra $<x>$ of $A$
generated by any element $x\in A$ is associative. $\ $Each alternative
algebra is$\ $a\ flexible algebra and a power associative algebra.

Algebras $A_{t}$ of dimension $2^{t}\ $obtained by the Cayley-Dickson
process, described above, are central-simple, flexible and\textit{\ }power
associative for all $t\geq 1$ and, in general, are not division algebras for
all $t\geq 1$.\ But there are fields (for example, the rational function
field) on which, if we apply the Cayley-Dickson process, the resulting
algebras $A_{t}\ $are division \ algebras for all $t\geq 1.$ (See [Br; 67]
and [Fl; 13] ).

\begin{equation*}
\end{equation*}%
\begin{equation*}
\end{equation*}

\textbf{2. About} \textbf{Fibonacci Quaternions }

\begin{equation*}
\end{equation*}%
\begin{equation*}
\end{equation*}

Since the above described algebras are usually without division, finding
quickly examples of invertible elements in an arbitrary algebra obtained by
the Cayley-Dickson process appear to be a not easy problem. A partial
solution for \ generalized real Quaternion algebras can be found using
Fibonacci quaternions.

Let $\mathbb{H}\left( \alpha _{1},\alpha _{2}\right) $ be the generalized
real\ quaternion algebra. In this algebra, every element has the form $%
a=a_{1}\cdot 1+a_{2}e_{2}+a_{3}e_{3}+a_{4}e_{4},$ where $a_{i}\in \mathbb{R}%
,i\in \{1,2,3,4\}$.

\qquad\ \qquad

In [Ho; 61], the Fibonacci quaternions were defined to be the quaternions on
the form%
\begin{equation}
F_{n}=f_{n}\cdot 1+f_{n+1}e_{2}+f_{n+2}e_{3}+f_{n+3}e_{4},  \tag{2.1}
\end{equation}%
called the \ $n$th Fibonacci quaternions, where 
\begin{equation}
f_{n}=f_{n-1}+f_{n-2,}\ n\geq 2,\   \tag{2.2}
\end{equation}%
with $f_{0}=0,f_{1}=1,$ are Fibonacci numbers$.$

The norm formula \ for the \ $n$th Fibonacci quaternions is:

\begin{equation}
\boldsymbol{n}\left( F_{n}\right) =F_{n}\overline{F}_{n}=3f_{2n+3}, 
\tag{2.3}
\end{equation}%
where \ $\overline{F}_{n}=f_{n}\cdot
1-f_{n+1}e_{2}-f_{n+2}e_{3}-f_{n+3}e_{4} $ is the conjugate of the $F_{n}$
in the algebra $\mathbb{H}\left( \alpha _{1},\alpha _{2}\right) $ (see [Ho;
61])$.$ There are many authors which studied Fibonacci quaternions in the
real division quaternion algebra giving more and surprising new properties
(see [Sw; 73], [Sa-Mu; 82] and [Ha; 12], [Fl, Sh; 13], [Fl, Sh; 13(1)]).

$\medskip ~\ $

\textbf{Theorem 2.1.} ([Fl, Sh; 13] Theorem 2.4. ) \textit{The norm of the} $%
n$th \textit{Fibonacci quaternion} $F_{n}$ \textit{in a generalized
quaternion algebra is}%
\begin{equation}
\boldsymbol{n}\left( F_{n}\right) \text{=}h_{2n+2}^{1+2\alpha _{2},3\alpha
_{2}}\text{+}(\alpha _{1}\text{-}1)h_{2n+3}^{1+2\alpha _{2},\alpha _{2}}%
\text{-}2(\alpha _{1}\text{-}1)\left( 1\text{+}\alpha _{2}\right)
f_{n}f_{n+1}.\medskip  \tag{2.4}
\end{equation}

We know that the expression for the $n$th term of a Fibonacci element is 
\begin{equation}
f_{n}=\frac{1}{\sqrt{5}}[a^{n}-b^{n}]=\frac{a^{n}}{\sqrt{5}}[1-\frac{b^{n}}{%
\alpha ^{n}}],  \tag{2.5}
\end{equation}%
where $a=\frac{1+\sqrt{5}}{2}$ and $b=\frac{1-\sqrt{5}}{2}.$

Using the \ above notations, we can compute the following limit 
\begin{equation*}
\underset{n\rightarrow \infty }{\lim }\boldsymbol{n}\left( F_{n}\right) =%
\underset{n\rightarrow \infty }{\lim }(f_{n}^{2}+\alpha
_{1}f_{n+1}^{2}+\alpha _{2}f_{n+2}^{2}+\alpha _{1}\alpha _{2}f_{n+3}^{2})=
\end{equation*}%
\begin{equation*}
\underset{}{=\underset{n\rightarrow \infty }{\lim }(\frac{a^{2n}}{5}\text{+}%
\alpha _{1}\frac{a^{2n+2}}{5}\text{+}\alpha _{2}\frac{a^{2n+4}}{5}\text{+}%
\alpha _{1}\alpha _{2}\frac{a^{2n+6}}{5})=}
\end{equation*}%
\begin{equation*}
=sgnE(\alpha _{1},\alpha _{2})\cdot \infty .
\end{equation*}%
Since $a^{2}=a+1,$ we $\ $have $E(\alpha _{1},\alpha _{2})=(\frac{1}{5}+%
\frac{\alpha _{1}}{5}a^{2}+\frac{\alpha _{2}}{5}a^{4}+\frac{\alpha
_{1}\alpha _{2}}{5}a^{6})=$\newline
$=\frac{1}{5}\left( 1+\alpha _{1}\left( a+1\right) +\alpha _{2}\left(
3a+2\right) +\alpha _{1}\alpha _{2}\left( 8a+5\right) \right) =$\newline
$=\frac{1}{5}[1+\alpha _{1}+2\alpha _{2}+5\alpha _{1}\alpha _{2}+a\left(
\alpha _{1}+3\alpha _{2}+8\alpha _{1}\alpha _{2}\right) ].$

If \ $E(\alpha _{1},\alpha _{2})>0,$ there exist a number $n_{1}\in \mathbb{N%
}$ such that for all\newline
$n\geq n_{1}$ we have 
\begin{equation*}
h_{2n+2}^{1+2\alpha _{2},3\alpha _{2}}+(\alpha _{1}-1)h_{2n+3}^{1+2\alpha
_{2},\alpha _{2}}-2(\alpha _{1}-1)\left( 1+\alpha _{2}\right) f_{n}f_{n+1}>0.
\end{equation*}%
If $E(\alpha _{1},\alpha _{2})<0,$ there exist a number $n_{2}\in \mathbb{N}$
such that for all $n\geq n_{2}$ we have 
\begin{equation*}
h_{2n+2}^{1+2\alpha _{2},3\alpha _{2}}+(\alpha _{1}-1)h_{2n+3}^{1+2\alpha
_{2},\alpha _{2}}-2(\alpha _{1}-1)\left( 1+\alpha _{2}\right) f_{n}f_{n+1}<0.
\end{equation*}

It results that for \ all $\alpha _{1},\alpha _{2}\in \mathbb{R}$ with $%
E(\alpha _{1},\alpha _{2})\neq 0,$ in the algebra $\mathbb{H}\left( \alpha
_{1},\alpha _{2}\right) $ there is a natural number $n_{0}=\max
\{n_{1},n_{2}\}$ such that $\boldsymbol{n}\left( F_{n}\right) \neq 0,$ hence 
$F_{n}$ is an invertible element, for all $n\geq n_{0}.$

In this way, Fibonacci Quaternion elements can provide us many important
information in the algebra $\mathbb{H}\left( \alpha _{1},\alpha _{2}\right) $
providing sets of invertible elements in algebraic structures without
division. For other details, see [Fl, Sh; 13].

\begin{equation*}
\end{equation*}%
\begin{equation*}
\end{equation*}

\textbf{3.} \textbf{Multiplication table in Cayley-Dickson algebras}%
\begin{equation*}
\end{equation*}%
\begin{equation*}
\end{equation*}

Multiplication table for algebras obtained by the Cayley-Dickson process
over the real field was studied in [Ba; 09]. In this paper, the author gave
an algorithm to find quickly product of two elements in these algebras. In
the following, we shortly present this algorithm. In [Ba; 13], the author
gave all $32$ possibilities to define a "Cayley-Dickson product" used in the
Cayley-Dickson doubling process, such that the obtained algebras are
isomorphic.

If we consider multiplication $\left( 1.1\right) $ under the form 
\begin{equation}
A\oplus A:\left( a_{1},a_{2}\right) \left( b_{1},b_{2}\right) :=\left(
a_{1}b_{1}+\gamma b_{2}\overline{a_{2}},\overline{a_{1}}b_{2}+b_{1}a_{2}%
\right) ,  \tag{3.1}
\end{equation}%
the obtained algebras are isomorphic with those obtained with multiplication 
$\left( 1.1\right) .\medskip $

For $\alpha _{1}=...=\alpha _{t}=-1$ and $K=\mathbb{R},$ in [Ba; 09] the
author described \ how we can multiply the basis vectors in the algebra $%
A_{t},\dim A_{t}=2^{t}=n$. He used the binary decomposition for the
subscript indices.

Let $\ e_{p},e_{q}$ be two vectors in the basis $B$ with $p,q$ representing
the binary decomposition for the indices of the vectors, that means $p,q$
are in $\mathbb{Z}_{2}^{n}.$ We have that $e_{p}e_{q}=\gamma _{n}\left(
p,q\right) e_{p\otimes q},$ where:

i) $p\otimes q$ are the "\textit{exclusive or}" for the binary numbers $p$
and $q$ (the sum of \ $p$ and $q$ in the group $\mathbb{Z}_{2}^{n});$

ii) $\gamma _{n}:\mathbb{Z}_{2}^{n}\times \mathbb{Z}_{2}^{n}\rightarrow
\{-1,1\}$ is a map, called the \textit{twist map}$.$

In this section, we will consider $K=\mathbb{R}.$ Using the same notations
as in the Bales's paper, we consider the following matrices: 
\begin{equation}
A_{0}=A=\left( 
\begin{array}{cc}
1 & 1 \\ 
1 & -1%
\end{array}%
\right) ,\quad B=\left( 
\begin{array}{cc}
1 & -1 \\ 
1 & 1%
\end{array}%
\right) ,\quad C=\left( 
\begin{array}{cc}
1 & -1 \\ 
-1 & -1%
\end{array}%
\right) .  \tag{3.2}
\end{equation}

In the same paper [Ba; 09], the author find the properties of the twist map $%
\gamma _{n}$ and put the signs of this map in a table. He partitioned the
twist table for \ $\mathbb{Z}_{2}^{n}$ into $2\times 2$ matrices and
obtained the following result:\medskip 
\begin{equation*}
\end{equation*}%
\begin{equation*}
\end{equation*}

\textbf{Theorem 3.1. ([}Ba; 09\textbf{]}, Theorem 2.2., p. 88-91\textbf{)} 
\textit{For} $n>0,$ \textit{the Cayley-Dickson twist table} $\gamma _{n}$ 
\textit{can be partitioned in quadratic} \textit{matrices of dimension} $2$ 
\textit{of the form} $A,B,C,-B,-C,$ \textit{defined in the relation} (3.2).

\begin{figure}[htbp]
\begin{center}
\includegraphics[height=5cm ,width=3cm, angle=270]{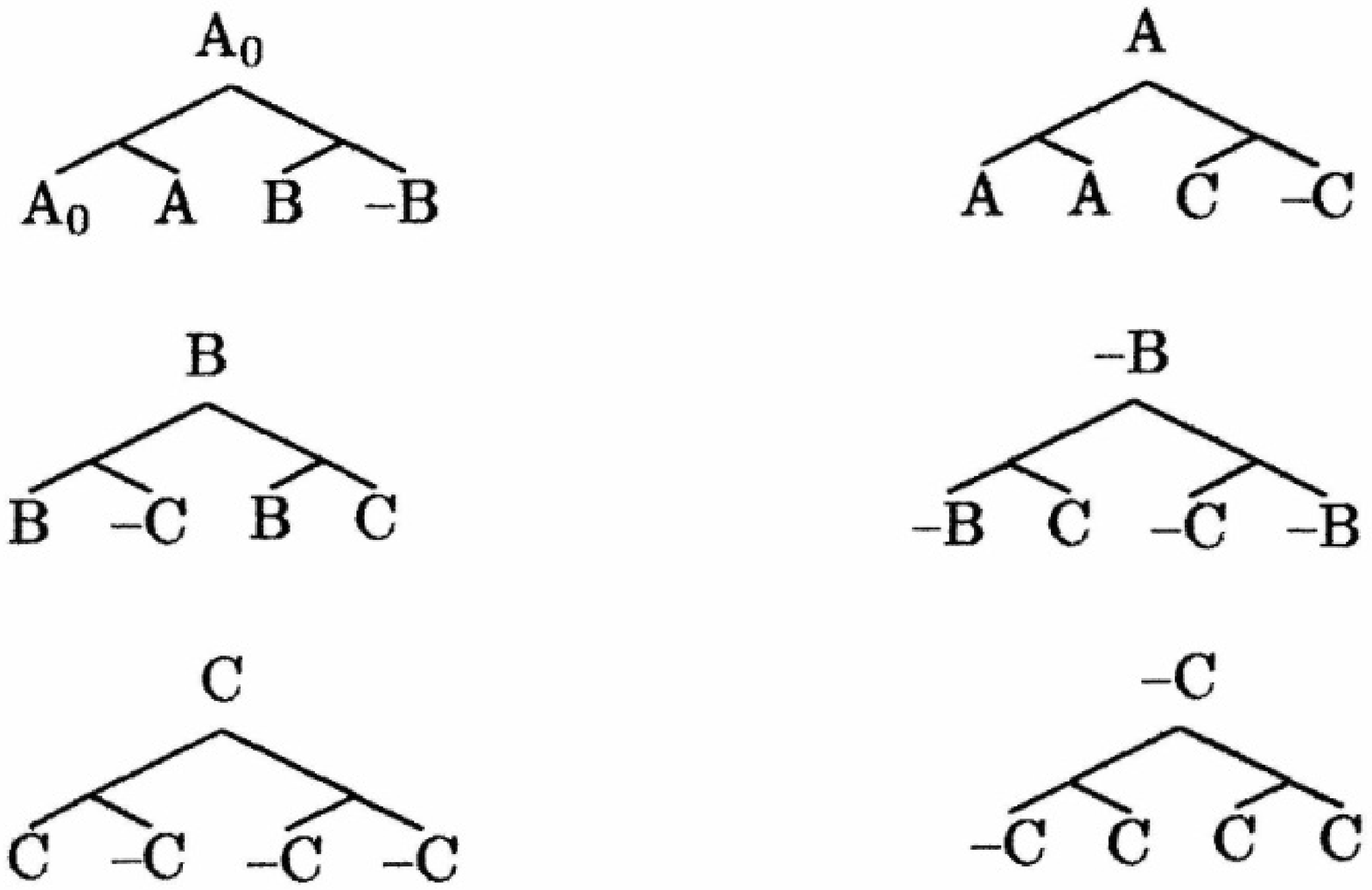}
\par
Fig. 1: Twist trees([Ba; 09], Table 9)
\end{center}
\end{figure}

\bigskip

\textbf{Definition 3.2.} Let $x=x_{0},x_{1},x_{2},....$ and~ $%
y=y_{0},y_{1},y_{2},.....$ be two sequences of real numbers. The ordered
pair 
\begin{equation*}
\left( x,y\right) =x_{0},y_{0},x_{1},y_{1},x_{2},y_{2},....
\end{equation*}%
is a sequence obtained by \textit{shuffling} \ the sequences $x$ and$~\
y.\medskip $

\textbf{Proposition 3.3.} \textit{Let} $A_{t}=\left( \frac{-1,...,-1}{%
\mathbb{R}}\right) $ \textit{be an algebra obtained by the Cayley-Dickson
process with multiplication given by relation }$\left( 3.1\right) ~$\textit{%
and} $\{e_{0}=1,e_{1},...,e_{n-1}\},\,$\newline
$\,n=2^{t}$ \textit{\ a basis in }$A_{t}$\textit{. Let} $\ r\geq
1,\,\,r<k\leq i<t.$ \textit{We have}%
\begin{equation}
\begin{tabular}{l|ll}
$\cdot $ & $e_{T}$ & $e_{T+1}$ \\ \hline
$e_{2^{k-r+1}}$ & $\left( -1\right) ^{r+2}e_{M}$ & $-\left( -1\right)
^{r+2}e_{M+1}$ \\ 
$e_{2^{k-r+1}+1}$ & $-\left( -1\right) ^{r+2}e_{M+1}$ & $-\left( -1\right)
^{r+2}e_{M}$%
\end{tabular}%
,  \tag{3.3}
\end{equation}%
\textit{where} \textit{the binary decomposition of} $M$ \textit{is} $%
M_{2}=2^{k}\otimes T.\medskip $

\textbf{Proof. }We \ compute $e_{2^{k-r+1}}e_{T}.$ We have $%
e_{2^{k-r+1}}e_{T}=\gamma \left( s,q\right) e_{M},$ where the binary
decomposition of $M$ is$\ M_{2}=2^{k-r+1}\otimes T$ and $s$ is the binary
decomposition for $2^{k-r+1}$ and $q$ is the binary decomposition for $T,$ 
\begin{equation*}
s=\underset{i-k+r-1}{\underbrace{00...0}}\underset{k-r+2}{\underbrace{100...0%
}},q=1\underset{i-k-1}{\underbrace{00...0}}\underset{k-r+1}{\underbrace{%
111..1}}\underset{r}{\underbrace{0...0}}.
\end{equation*}
By "shuffling" $s\otimes q$, it results 
\begin{equation*}
\underset{i-k}{\underbrace{01~00~00...00}\ }\underset{k-2r-1}{\underbrace{%
01~01~01~...01~}}\underset{r+2}{\ \underbrace{11~01~01~...01}}\underset{r}{\ 
\underbrace{00~00~...00~00}}.
\end{equation*}

Starting with $A_{0},$we get:

\begin{equation*}
\underset{i-k}{\underbrace{A_{0}\overset{01}{\rightarrow }A\overset{00}{%
\rightarrow }\text{...}\overset{00}{\rightarrow }}}\underset{k-2r-1}{%
\underbrace{A\overset{01}{\rightarrow }A\overset{01}{\rightarrow }\text{...}%
\overset{01}{\rightarrow }}}\underset{r+2}{\underbrace{A\overset{11}{%
\rightarrow }-C\overset{01}{\rightarrow }C\overset{01}{\rightarrow }-C%
\overset{01}{\rightarrow }C\text{...}\overset{01}{\rightarrow }\left(
-1\right) ^{r+2}C}}\underset{r}{\underbrace{\overset{00}{\rightarrow }\text{%
...}\overset{00}{\rightarrow }\left( -1\right) ^{r+2}C}}.
\end{equation*}%
Therefore $\gamma \left( s,q\right) =\left( -1\right) ^{k-r+1}.$

Now, we compute $e_{2^{k-r+1}}e_{T+1}.$ For this, we will "shuffling" $%
\underset{i-k+r-1}{\underbrace{00...0}}\underset{k-r+2}{\underbrace{100...0}}
$ with $1\underset{i-k-1}{\underbrace{00...0}}\underset{k-r+1}{\underbrace{%
111..1}}\underset{r}{\underbrace{0...1}}.$ It results 
\begin{equation*}
\underset{i-k}{\underbrace{01~00~00...00}\ }\underset{k-2r-1}{\underbrace{%
01~01~01~...01~}}\underset{r+2}{\ \underbrace{11~01~01~...01}}\underset{r}{\ 
\underbrace{00~00~...00~01}}.
\end{equation*}
Starting with $A_{0},$we get:

\begin{equation*}
\underset{i-k}{\underbrace{A_{0}\overset{01}{\rightarrow }A\overset{00}{%
\rightarrow }\text{...}\overset{00}{\rightarrow }}}\underset{k-2r-1}{%
\underbrace{A\overset{01}{\rightarrow }A\overset{01}{\rightarrow }\text{...}%
\overset{01}{\rightarrow }}}\underset{r+2}{\underbrace{A\overset{11}{%
\rightarrow }-C\overset{01}{\rightarrow }C\overset{01}{\rightarrow }-C%
\overset{01}{\rightarrow }C\text{...}\overset{01}{\rightarrow }\left(
-1\right) ^{r+2}C}}\underset{r}{\underbrace{\overset{00}{\rightarrow }\text{%
...}\overset{01}{\rightarrow }\left( -1\right) ^{r+3}C}}.
\end{equation*}

For $e_{2^{k-r+1}+1}e_{T},$ "shuffling" $\underset{i-k+r-1}{\underbrace{%
00...0}}\underset{k-r+2}{\underbrace{100...1}}$ with $1\underset{i-k-1}{%
\underbrace{00...0}}\underset{k-r+1}{\underbrace{111..1}}\underset{r}{%
\underbrace{0...0}},$ it results 
\begin{equation*}
\underset{i-k}{\underbrace{01~00~00...00}\ }\underset{k-2r-1}{\underbrace{%
01~01~01~...01~}}\underset{r+2}{\ \underbrace{11~01~01~...01}}\underset{r}{\ 
\underbrace{00~00~...00~10}}.
\end{equation*}%
$\ $Starting with $A_{0},$we get: 
\begin{equation*}
\underset{i-k}{\underbrace{A_{0}\overset{01}{\rightarrow }A\overset{00}{%
\rightarrow }\text{...}\overset{00}{\rightarrow }}}\underset{k-2r-1}{%
\underbrace{A\overset{01}{\rightarrow }A\overset{01}{\rightarrow }\text{...}%
\overset{01}{\rightarrow }}}\underset{r+2}{\underbrace{A\overset{11}{%
\rightarrow }-C\overset{01}{\rightarrow }C\overset{01}{\rightarrow }-C%
\overset{01}{\rightarrow }C\text{...}\overset{01}{\rightarrow }\left(
-1\right) ^{r+2}C}}\underset{r}{\underbrace{\overset{00}{\rightarrow }\text{%
...}\overset{10}{\rightarrow }\left( -1\right) ^{r+3}C}}.
\end{equation*}%
For $e_{2^{k-r+1}+1}e_{T+1},$ we compute first $\left( 2^{k-r+1}+1\right)
\otimes \left( T+1\right) .$ We obtain:\newline
\begin{equation*}
\left( 2^{k-r+1}+1\right) \otimes \left( T+1\right) =\newline
\end{equation*}%
\begin{equation*}
=\left( \underset{i-k+r-1}{\underbrace{00...0}}\underset{k-r+2}{\underbrace{%
100...1}}\right) \otimes \left( 1\underset{i-k-1}{\underbrace{00...0}}%
\underset{k-r+1}{\underbrace{111..1}}\underset{r}{\underbrace{0...1}}\right)
=\medskip \medskip \newline
\end{equation*}%
\begin{equation*}
=\ \ \ \underset{i-k}{\underbrace{10...0}}\underset{r-1}{\underbrace{11..1}}%
0~\underset{k-2r+1}{\underbrace{1...1}}\underset{r}{\underbrace{0...0}}\
=2^{k-r+1}\otimes T=M.
\end{equation*}%
$\medskip $

Now, "shuffling" $\underset{i-k+r-1}{\underbrace{00...0}}\underset{k-r+2}{%
\underbrace{100...1}}$ with $1\underset{i-k-1}{\underbrace{00...0}}\underset{%
k-r+1}{\underbrace{111..1}}\underset{r}{\underbrace{0...1}},$ it results 
\begin{equation*}
\underset{i-k}{\underbrace{01~00~00...00}\ }\underset{k-2r-1}{\underbrace{%
01~01~01~...01~}}\underset{r+2}{\ \underbrace{11~01~01~...01}}\underset{r}{\ 
\underbrace{00~00~...00~11}}
\end{equation*}%
\newline
Starting with $A_{0},$we get:

\begin{equation*}
\underset{i-k}{\underbrace{A_{0}\overset{01}{\rightarrow }A\overset{00}{%
\rightarrow }\text{...}\overset{00}{\rightarrow }}}\underset{k-2r-1}{%
\underbrace{A\overset{01}{\rightarrow }A\overset{01}{\rightarrow }\text{...}%
\overset{01}{\rightarrow }}}\underset{r+2}{\underbrace{A\overset{11}{%
\rightarrow }-C\overset{01}{\rightarrow }C\overset{01}{\rightarrow }-C%
\overset{01}{\rightarrow }C\text{...}\overset{01}{\rightarrow }\left(
-1\right) ^{r+2}C}}\underset{r}{\underbrace{\overset{00}{\rightarrow }\text{%
...}\overset{11}{\rightarrow }\left( -1\right) ^{r+3}C}}.
\end{equation*}%
$\Box \medskip $%
\begin{equation*}
\end{equation*}

\textbf{4. Some applications in Algebra and Coding Theory}%
\begin{equation*}
\end{equation*}

Let $A_{t}$ be an algebra obtained by the Cayley-Dickso process over the
field $\mathbb{R}$, with the basis $\{1,e_{2},...,e_{n}\},n=2^{t}.$ The unit
elements in $A_{t}$ are $\{\pm 1,\pm e_{2},...,\pm e_{n}\}.$ In [Ma, Be, Ga;
09], the authors defined \textit{the integers} of the $A_{t}$ to be the set 
\begin{equation*}
A_{t}[\mathbb{Z}]=\{x_{1}\cdot 1+\underset{i=2}{\overset{2^{n}}{\sum }}%
x_{i}\cdot e_{i},\ x_{1},x_{i}\in \mathbb{Z},i\in \{2,...,n\}\}.
\end{equation*}

\begin{equation*}
\end{equation*}%
$A_{t}[\mathbb{Z}]$ is a non-associative and non-commutative ring on which
the following equivalence relation can be defined.\medskip

\textbf{Definition 4.1.} Let $a,x,y\in A_{t}[\mathbb{Z}].$ We say that $x,y$ 
\textit{are right(left) congruent modulo} $a$ if and only if there is the
element $b\in A_{t}[\mathbb{Z}]$ such that 
\begin{equation}
x-y=ba~(\text{ or }x-y=ab).  \tag{4.1}
\end{equation}

We denote this relation with $x\equiv _{r}y$ \textit{mod} $a$ (or $x\equiv
_{s}y$ \textit{mod} $a)$ and this relation is well defined. We will consider
the quotient ring 
\begin{equation*}
A_{t}[\mathbb{Z}]_{a}=\{x\text{ \textit{mod~}}a~/~x\in A_{t}[\mathbb{Z}]\}.
\end{equation*}%
If $a\neq 0$ is not a zero divisor, then $A_{t}[\mathbb{Z}]_{a}$ has $%
N\left( a\right) ^{2^{n-1}}$ elements (see [Ma, Be, Ga; 09] \ for other
details).

Since algebras $A_{t}$ are poor in properties, due to the
power-associativity, if we take $w\in A_{t}[\mathbb{Z}]$, then the set $%
\mathbb{U}=\{a+bw~/~a,b\in \mathbb{Z}\}$ become an associative and a
commutative ring with $\mathbb{U}\subset A_{t}[\mathbb{Z}].$

Let $\mathbb{U}$ be the ring defined above, included in $A_{t}[\mathbb{Z}],$ with $%
t\in \{2,3\}.\medskip $

\textbf{Definition 4.2.} An element $x\in \mathbb{U}$ \textit{is prime} in $%
\mathbb{U}$ if $x$ is not an invertible element in $\mathbb{U}$ and if $x=ab$%
, then $a$ or $b$ is an invertible element in $\mathbb{U}$.\medskip

It is obvious that if $\pi \in \mathbb{U}$ is a prime element, then $n\left(
\pi \right) $ is a prime element in $\mathbb{Z}.$

If we consider relation $\left( 4.1\right) \,$\ on $\mathbb{U},$ due to
commutativity, "the left" is the same with "the right" and if $\pi $ is a
prime element in $\mathbb{U},$ therefore $\mathbb{U}_{\pi }$ is a field
isomorphic with $\mathbb{Z}_{p},$ where $n\left( \pi \right) =p,p$ a prime
element in $\mathbb{Z},$ as we can see from the above statements.\medskip

\textbf{Proposition 4.3. (}[Fl; 14], [Gu; 13], [Hu; 94]\textbf{\ )}

\textit{i)} \textit{If} $x,y\in \mathbb{V}$\textit{, then there are} $z,v\in 
\mathbb{V}$ \textit{such that} $x=zy+v,$ \textit{with} $N\left( v\right)
<N\left( y\right) .$

\textit{ii) With the above notation, we have that the remainder} $v$ \textit{%
has the formula}%
\begin{equation}
v=x-\left[ \frac{x\overline{y}}{y\overline{y}}\right] y,  \tag{4.2}
\end{equation}%
\textit{where the symbol} $[~,]$ \textit{is the rounding to the closest
integer. For the octonions, the rounding of an octonion integer can be found
by rounding the coefficients of the basis, separately, to the closest
integer.}\medskip

\textbf{Proposition 4.4. }([Fl; 14], [Gu; 13], [Hu; 94])

\textit{i) The above relation is an equivalence relation on }$\mathbb{U}$%
\textit{. The set of equivalence class is denoted by} $\mathbb{U}_{\pi }$ 
\textit{and is called\ the residue class of } $\mathbb{U}$\textit{\ modulo} $%
\pi .$

\textit{ii) The modulo function }$\mu :\mathbb{U\rightarrow U}_{\pi }$ 
\textit{is} $\mu \left( x\right) =v\ $mod $\pi =x-\left[ \frac{x\overline{y}%
}{y\overline{y}}\right] y,$ where $x=z\pi +v,$ \textit{with} $N\left( \pi
\right) <N\left( y\right) .$

\textit{iii)} $\mathbb{U}_{\pi }$ \textit{is a field } \textit{isomorphic
with} $\mathbb{Z}_{p},p=N(\pi )\medskip ,p$ \textit{a prime number.\medskip }

\textbf{Remark 4.6. ([}Ne, In,Fa, El, Pa; 01\textbf{])} From the above, we
have that for$~v_{i},v_{j}\in \mathbb{U}_{\pi },i,j\in \{1,2,...,p-1\},$ $%
u_{i}+u_{j}=u_{k}~$if and only if $k=i+j$ \textit{mod} $p$ and $u_{i}\cdot
u_{j}=u_{k}~$if and only if $k=i\cdot j$ \textit{mod} $p.$ From here, we
have the following labelling procedure:

1) Let $\pi $ $\in \mathbb{U}$ be a prime, with $n\left( \pi \right) =p,p$ a
prime number, $\pi =a+bw,a,b\in \mathbb{Z}.$

2) Let $s\in \mathbb{Z}$ be the only solution of the equation $a+bx$ \textit{%
mod} $p,~x\in \{0,1,2,...,p-1\}.$

3) The element $k\in \mathbb{Z}_{p}$ is the label of the element $u=m+nw\in 
\mathbb{U}$ if $m+ns=k$ \textit{mod} $p$ and $n\left( u\right) $ is minimum.

In this way, we obtain the map%
\begin{equation*}
\alpha :\mathbb{Z}_{p}\rightarrow \mathbb{U}_{\pi },~\alpha \left( \mathbf{m}%
\right) =\mu \left( m+\pi \right) =\left( m+\pi \right) ~\text{\textit{mod }}%
\pi .
\end{equation*}

\textbf{Example 4.5.}

Let $t=2,w=1+e_{2}+e_{3}+e_{4},p=13,\pi =-1+2w.$ We remark that $n\left( \pi
\right) =13$ and $w^{2}-2w+4=0.$ The field $\mathbb{U}_{\pi }$ isosmorphic
with $\mathbb{Z}_{13}$ is 
\begin{equation*}
\mathbb{U}_{\pi }=\{0,1,2,3,-3+w,-2+w,-1+w,1-w,2-w,3-w,-3,-2,-1\}
\end{equation*}%
Indeed, using relations $w^{2}=2w-4$ and $\overline{w}=2-w,$ we have:

$\mathbf{4}=4+\pi =3+2w=-3+w,$ since $3+2w=\left( -1+2w\right) \overline{w}%
+w-3,$ with $n\left( w-3\right) =7<13=n\left( \pi \right) ;$

$\mathbf{5}=5+\pi =4+2w=-2+w;$

$\mathbf{6}=6+\pi =5+2w=-1+w.$

Using the above labelling procedure, we have \newline
$\alpha :\mathbb{Z}_{p}\rightarrow \mathbb{U}_{\pi },\alpha \left( 0\right)
=0,\alpha \left( 1\right) =1,\alpha \left( 2\right) =2,\alpha \left(
3\right) =3,$\newline
$\alpha \left( 4\right) =-3+w,\alpha \left( 5\right) =-2+w,\alpha \left(
6\right) =-1+w,\alpha \left( 7\right) =1-w,$\newline
$\alpha \left( 8\right) =2-w.\alpha \left( 9\right) =3-w,\alpha \left(
10\right) =-3,\alpha \left( 11\right) =-2,\alpha \left( 12\right)
=-1.\medskip $

\textbf{Remark 4.6.} Since each natural number can be write as a sum of four
squares, if $m\in \mathbb{N},$ such that $%
m=a_{1}^{2}+a_{2}^{2}+a_{3}^{2}+a_{4}^{2},a_{i}\in \mathbb{N},i\in \{1,2,3\}$
and if $q=2a_{1},$ therefore the equation 
\begin{equation}
x^{2}-qx+m=0,  \tag{4.3}
\end{equation}%
has always solutions in $A_{t},$ for all $t.$ Indeed, let $z=a_{1}\cdot
1+a_{2}\cdot e_{i}+a_{3}\cdot e_{j}+a_{4}\cdot e_{k},$ where $i\neq j\neq k$
and $e_{i},e_{j},e_{k}\in \{e_{2},...,e_{n}\},n=2^{t}.$ The element $z$ is
always a solution of the equation $(4.3),$ since $t\left( z\right) =2a_{1}=q$
and $n\left( x\right) =a_{1}^{2}+a_{2}^{2}+a_{3}^{2}+a_{4}^{2}=m.\medskip $

\textbf{Remark 4.7.} Such a field obtained above have many applications in
Coding Theory, since on these fields can be constructed good codes which can
detected and corrected some error patterns which occur most frequently (see
[Fl; 14], [Gu; 13], [Hu; 94], [Ma, Be, Ga; 09], [Ne, In,Fa, El, Pa; 01]).%
\begin{equation*}
\end{equation*}

\begin{equation*}
\end{equation*}

\textbf{References}%
\begin{equation*}
\end{equation*}

[Al; 98] S.M. Alamouti, \textit{A simple transmit diversity technique for
wireless communications}, IEEE \ J. Selected Areas Communications, \textbf{%
16(1998)}, 1451-1458.

[Ba; 13] J.W. Bales, \textit{A Catalog of \ Cayley-Dickson-like Products}, 
\newline
http://arxiv.org/pdf/1107.1301v4.pdf.

[Ba; 09] J. W. Bales, \ \textit{A Tree for Computing the Cayley-Dickson Twist%
}, \ Missouri J. Math. Sci., \textbf{21(2)(2009)}, 83--93.

[Br; 67] R. B. Brown, \textit{On generalized Cayley-Dickson algebras},
Pacific J. of Math.,\textbf{\ 20(3)(1967)}, 415-422.

[Fl; 14] C. Flaut, \textit{Codes over a subset of Octonion Integers}, 
\newline
http://arxiv.org/pdf/1401.7828.pdf

[Fl; 13] C. Flaut, \textit{Levels and sublevels of \ algebras obtained by
the Cayley-Dickson process}, Ann Mat Pur Appl \textbf{192(2013)}, 1099-1114.

[Fl, Sh; 13] C. Flaut, V. Shpakivskyi, \ \textit{On Generalized Fibonacci
Quaternions and Fibonacci-Narayana Quaternions}, Adv. Appl. Cliff ord
Algebras \textbf{23(2013)}, 673--688.

[Fl, Sh; 13(1)] C. Flaut, V. Shpakivskyi, \ \textit{Real Matrix
Representations for the Complex Quaternions}, Adv. Appl. Cliff ord Algebras 
\textbf{23(2013)}, 657--671.

[Gu; 13] M. G\"{u}zeltepe, \textit{Codes over Hurwitz integers},
DiscreteMath, \textbf{313(5)(2013)}, 704-714.

[Ha; 12] S. Halici, \textit{On Fibonacci Quaternions}, Adv. in Appl.
Clifford Algebras, \textbf{22(2)(2012)}, 321-327.

[Ho; 61] A. F. Horadam, \textit{A Generalized Fibonacci Sequence}, Amer.
Math. Monthly, \textbf{68(1961)}, 455-459.

[Hu; 94] K. Huber, \textit{Codes over Gaussian integers}, IEEE Trans.
Inform. Theory,\textbf{\ 40(1994)}, 207--216.

[Ma, Be, Ga; 09] C. Martinez, R. Beivide, E. Gabidulin, \textit{Perfect
codes from Cayley graphs over Lipschitz integers}, IEEE Trans. Inform.
Theory \textbf{55(8)(2009),} 3552--3562.

[Ne, In,Fa, El, Pa; 01] T.P. da N. Neto, J.C. Interlando, M.O. Favareto, M.
Elia, R. Palazzo Jr., \textit{Lattice constellation and codes from quadratic
number fields}, IEEE Trans. Inform.Theory \textbf{47(4)(2001)} 1514--1527.

[Sa-Mu; 82] P. V. Satyanarayana Murthy, \textit{Fibonacci-Cayley Numbers},
The Fibonacci Quarterly, \textbf{20(1)(1982)}, 59-64.

[Sc; 66] R. D. Schafer, \textit{An Introduction to Nonassociative Algebras,}
Academic Press, New-York, 1966.

[Sc; 54] R. D. Schafer, \ \textit{On the algebras formed by the
Cayley-Dickson process,} Amer. J. Math., \textbf{76(1954)}, 435-446.

[Sw; 73] M. N. S. Swamy, \textit{On generalized Fibonacci Quaternions}, The
Fibonacci Quaterly, \textbf{11(5)(1973)}, 547-549.%
\begin{equation*}
\end{equation*}

Cristina FLAUT

{\small Faculty of Mathematics and Computer Science,}

{\small Ovidius University,}

{\small Bd. Mamaia 124, 900527, CONSTANTA,}

{\small ROMANIA}

{\small http://cristinaflaut.wikispaces.com/}

{\small http://www.univ-ovidius.ro/math/}

{\small e-mail:}

{\small cflaut@univ-ovidius.ro}

{\small cristina\_flaut@yahoo.com}

\end{document}